\documentclass[twoside, 10pt]{article}
\usepackage{amssymb, amsmath, mathrsfs, amsthm}
\usepackage{graphicx}
\usepackage{color}
\usepackage{float}

\input{epsf}

\newcommand{\bi}{\begin{itemize}}
\newcommand{\ei}{\end{itemize}}
\newcommand{\bt}{\begin{tabular}}
\newcommand{\et}{\end{tabular}}
\newcommand{\beq}{\begin{equation}}
\newcommand{\eeq}{\end{equation}}
\newcommand{\beqs}{\begin{eqnarray*}}
\newcommand{\eeqs}{\end{eqnarray*}}
\newcommand{\be}{\begin{enumerate}}
\newcommand{\ee}{\end{enumerate}}

\definecolor{DarkGreen}{rgb}{0.2, 0.6, 0.3}

\newcommand{\Gk}{G$_k$}

\newtheorem{theorem}{Theorem}

\newtheorem{proposition}{Proposition}

\begin{document}
\title{\bf Note on a generalization of Gallai-Ramsey numbers}
\author{Colton Magnant\footnote{Department of Mathematics, Clayton State University, Morrow, GA, 30260, USA.} \footnote{Academy of Plateau Science and Sustainability, Xining, Qinghai 810008, China}, Zhuojun Magnant\footnote{Department of Mathematical Sciences, Georgia Southern University, Statesboro, GA 30460, USA.}}
\maketitle

\begin{abstract}
A colored complete graph is said to be Gallai-colored if it contains no rainbow triangle. This property has been shown to be equivalent to the existence of a partition of the vertices (of every induced subgraph) in which at most two colors appear on edges between the parts and at most one color appears on edges in between each pair of parts. We extend this notion by defining a coloring of a complete graph to be $k$-Gallai if every induced subgraph has a nontrivial partition of the vertices such that there are at most $k$ colors present in between parts of the partition.  The generalized $(k, \ell)$ Gallai-Ramsey number of a graph $H$ is then defined to be the minimum number of vertices $N$ such that every $k$-Gallai coloring of a complete graph $K_{n}$ with $n \geq N$ using at most $\ell$ colors contains a monochromatic copy of $H$.  We prove bounds on these generalized $(k, \ell)$ Gallai-Ramsey numbers based on the structure of $H$, extending recent results for Gallai colorings.
\end{abstract}

\section{Introduction}

Unless specifically stated, all colorings in this paper are colorings of the edges of complete graphs.  For notation not defined here, we refer the reader to \cite{CLZ11}.

A coloring of a complete graph $G$ is called a \emph{Gallai coloring} if it contains no rainbow triangle in honor of the following seminal theorem of Gallai \cite{G67}, also attributed to Cameron and Edmonds \cite{CE97}.  See \cite{GS04} for a restatement with additional comments.

\begin{theorem}[\cite{CE97, G67, GS04}]\label{Thm:G}
In every Gallai colored complete graph, there is a partition of the vertices such that there are at most two colors on the edges between the parts and only one color on edges between each pair of parts.
\end{theorem}

A partition resulting from Theorem~\ref{Thm:G} is called a \emph{Gallai partition}.  Given a positive integer $k$ and graphs $G$ and $H$, the \emph{$k$-color Gallai-Ramsey number of $G$ and $H$}, denoted $gr_{k}(G : H)$, is defined to be the minimum number of vertices $N$ such that every coloring of a complete graph on at least $N$ vertices using at most $k$ colors must contain either a rainbow copy of $G$ or a monochromatic copy of $H$.  Note that the study of Gallai colorings is important for the case when $G = K_{3}$.

Using Theorem~\ref{Thm:G}, Gy{\'a}rf{\'a}s et al. \cite{GSS08} were able to prove the following result determining the asymptotic behavior of Gallai-Ramsey numbers as a function of the number of colors.

\begin{theorem}[\cite{GSS08}]\label{Thm:G-Asymp}
If $H$ is bipartite and not a star, then $gr_{k}(K_{3} : H)$ grows as a linear function of $k$.  If $H$ is not bipartite, then $gr_{k}(K_{3} : H)$ grows as an exponential function of $k$.
\end{theorem}

A generalization of Theorem~\ref{Thm:G} was proven in \cite{FM12}.

\begin{theorem}[\cite{FM12}]\label{Thm:P}
Of every coloring of a complete graph containing no rainbow triangle with a pendant edge, there is a partition of the vertices such that either
\bi
\item there are at most two colors on the edges between the parts, or
\item there are at most three colors on the edges between the parts and only one color on edges between each pair of parts.
\ei
\end{theorem}

Although Theorem~\ref{Thm:P} provides a much weaker conclusion than Theorem~\ref{Thm:G}, Fujita and Magnant \cite{FM12} were able to produce an analog of Theorem~\ref{Thm:G-Asymp} for colorings free of rainbow copies of a triangle with a pendant edge.  For the sake of notation, let $S_{3}^{+}$ denote the triangle with a pendant edge.

\begin{theorem}[\cite{FM12}]\label{Thm:P-Asymp}
If $H$ is bipartite and not a star, then $gr_{k}(S_{3}^{+} : H)$ grows as a linear function of $k$.  If $H$ is not bipartite, then $gr_{k}(S_{3}^{+} : H)$ grows as an exponential function of $k$.
\end{theorem}

In this work, we assume a still weaker structure on the graph.  We say that a colored complete graph has a nontrivial \Gk-partition if there is a partition of the vertices into at least two parts such that there are at most $k$ colors on edges between the parts.  Say that a colored graph is \emph{$k$-Gallai} if every nontrivial induced subgraph has a nontrivial \Gk-partition.  We may then define the \emph{generalized $(k, \ell)$ Gallai-Ramsey number} of a graph $H$, denoted $ggr_{k, \ell}(H)$, to be the minimum number of vertices $N$ such that every $k$-Gallai coloring of a complete graph $K_{n}$, with $n \geq N$, using at most $\ell$ colors, contains a monochromatic copy of $H$.

Note that with $k = 2$, this structure is a bit stronger than the structure provided by Theorem~\ref{Thm:P} but still weaker than the structure provided by Theorem~\ref{Thm:G}.  In fact, Fujita and Magnant \cite{FM12} also proved the following more general result concerning a star $S_{t}$ with an additional edge between leaves.  We state the result here using our notation of $k$-Gallai colorings.

\begin{theorem}[\cite{FM12}]
Any rainbow $S_{k}^{+}$-free coloring of a complete graph is $3$-Gallai.
\end{theorem}

Within this general framework, we prove the following, our main result of this work.

\begin{theorem}\label{Thm:Main}
Suppose we are given positive integers $k$ and $\ell$ and a graph $H$.  If $H$ is bipartite and not a star, then $ggr_{k, \ell}(H)$ grows as a linear function of $\ell$.  If $H$ is not bipartite, then $ggr_{k, \ell}(H)$ grows as an exponential function of $\ell$, that grows with $\chi(H)$.
\end{theorem}

In the case when $H$ is a star, we have the following easy result when $\ell$ is assumed to be sufficiently large.  Although $\ell$ is assumed to be large, note that the number is not a function of $\ell$.

\begin{proposition}
Given integers $k, t \geq 2$, if $\ell$ is sufficiently large, then 
$$
ggr_{k,\ell} (S_t) = 2k(t-1)+1.
$$
\end{proposition}

\begin{proof}
For the lower bound, consider the graph on $2k(t-1)$ vertices consisting of two sets $A_1$ and $A_2$, each of order $k(t-1)$. Color the edges between $A_1$ and $A_2$ with $k$ colors so that each vertex is incident to exactly $t-1$ edges of each color (to the opposite set). Since $\ell$ is sufficiently large,  there exists a coloring within each subgraph induced on $A_i$ using the other $\ell-k$ colors (with the condition that each $A_i$ is $k$-Gallai) and avoiding a monochromatic copy of $S_t$. This graph is $k$-Gallai and has no monochromatic $S_t$.

For the upper bound, suppose $G$ is a $k$-Gallai coloring of $K_n$ using at most $\ell$ colors and suppose that $n \geq 2k(t-1)+1$. Consider a $G_k$-partition of $G$ and let $A$ be the smallest part of this partition. Since this is a nontrivial partition, $|A| \leq \lfloor{\frac{n}{2}} \rfloor \leq k(t-1)$ which means there are at least $k(t-1) + 1$ vertices in $G \setminus A$. If we let $v \in V(A)$, then $v$ has at most $k$ colors on its edges to $G \setminus A$. By the pigeonhole principle, one color must appear on at least $t$ edges from the vertex $v$ to $G \setminus A$. This forms a monochromatic $S_t$ centered at $v$ and completes the proof. \end{proof}

\section{Preliminaries}

Given integers $m, n \geq 2$, let $z(m; n)$ denote the minimum number of edges $q$ such that any balanced bipartite graph $B_{m, m}$ of order $2m$ containing at least $q$ edges contains a copy of $K_{n, n}$.  Bollob{\'a}s stated the following theorem in \cite{B04}. 

\begin{theorem}[\cite{Z63}]\label{Thm:Z-Numbers}  
$$
z(m; n) < (n - 1)^{1/n}m^{2 - 1/n} + \frac{n - 1}{2}m = o(m^{2}).
$$
\end{theorem}

Given positive integers $\ell, n_{1}, n_{2}, \dots, n_{\ell}$ with $\sum n_{i} = n$, let $L(t, n_{1}, n_{2}, \dots, n_{\ell})$ denote the coloring of $K_{n}$ constructed as follows.  Let $G_{1}$ be a coloring of $K_{n_{1}}$ entirely with color $1$.  Then for each $i \geq 1$, let $G_{i + 1}$ be a colored complete graph constructed by taking the graph $G_{i}$ with the addition of $n_{i + 1}$ vertices, each with all incident edges having color $i + 1$.  Then $G_{\ell}$ is the desired coloring and note that $G_{\ell}$ contains no rainbow triangle.  Further note that $G_{\ell}$ is $k$-Gallai for every $k \geq 1$.

\section{Proof of Theorem~\ref{Thm:Main}}

In this section, we prove Theorem~\ref{Thm:Main}.  If fact, we prove general bounds for each case separately.

First we consider the case when $H$ is bipartite and not a star.  Let $b$ denote the $k$-color bipartite Ramsey number of $H$, $b = b(k, H) = \min\{p : $ any $k$-coloring of $K_{p, p}$ contains a monochromatic copy of $H\}$.  Such numbers are known to exist by a result of Erd{\H o}s and Rado in \cite{ER56}.  %
Let $t = R_{k}(H) - 1$ and note that since $H$ is not a star, $t \geq 3$.

Let $n$ be the order of the larger part of the bipartite graph $H$.  By Theorem~\ref{Thm:Z-Numbers}, there is a number $z = \min\{M :$ any bipartite graph with partite sets of sizes $m_{1}, m_{2} \geq M$ containing at least $\frac{m_{1}m_{2}}{2k}$ %
edges contains a copy of $K_{n, n}\}$.  By the definition of $n$, such a monochromatic copy of $K_{n, n}$ would contain the desired monochromatic copy of $H$.

\begin{theorem}\label{Thm:MainBip}
Given positive integers $k \leq \ell$ and a bipartite graph $H$ that is not a star,
$$
n + \ell(m - 1) < ggr_{k, \ell}(H) \leq tb + (2k + 1)((z - 1)\ell + 1)(b - 1) + z.
$$
\end{theorem}

Since $t, b$ and $z$ are all functions of $k$ and $H$ but not $\ell$, this result shows that $ggr_{k, \ell}(H)$ is a linear function of $\ell$ when $H$ is bipartite.

\begin{proof}
Given the bipartite graph $H$, let $m$ be the order of the smaller part of $H$.  The lower bound is provided by $L(\ell, (n + m - 1), m - 1, m - 1, \dots, m - 1)$.  This graph is certainly $k$-Gallai, contains no monochromatic copy of $H$ as a subgraph, and has order $n + m - 1 + (\ell - 1)(m - 1) = n + \ell(m - 1)$.

Now we consider the upper bound.  For convenience, let $s = (z - 1)\ell$.  Let $G$ be an $\ell$-coloring of the complete graph on $tb + s(b - 1) + b$ vertices.  By assumption, there is a partition of $V(G)$ with at most $k$ colors on the edges between the parts.  Choosing one vertex from each part induces a colored complete graph using at most $k$ colors, so this implies that there are at most $t$ parts in the partition.  Let $H_{1}$ be the largest part in the partition and note that $|H_{1}| \geq \frac{|G|}{t} > b$.  If $|G \setminus H_{1}| \geq b$, then there exists a monochromatic copy of $H$ within the bipartite graph between $H_{1}$ and $G \setminus H_{1}$ (by the definition of $b$).  Hence, we may assume $|H_{1}| \geq |G| - b + 1$.

Let $v_{1}$ be a vertex in $G \setminus H_{1}$ and note that $v_{1}$ has at most $k$ colors on the edges to $H_{1}$.  Color $v_{1}$ with one of the colors that appears on at least $\frac{|H_{1}|}{k}$ edges to $H_{1}$.  Using $H_{1}$ in place of $G$, by assumption, there is a partition of $H_{1}$ and by the same argument, the largest part $H_{2}$ must have $|H_{2}| \geq |H_{1}| - b + 1$.  Let $v_{2}$ be a vertex of $H_{1} \setminus H_{2}$ and, as with $v_{1}$, color $v_{2}$ with one of the colors that appears on at least $\frac{|H_{2}|}{k}$ edges to $H_{2}$.

Repeat this process $s + 1$ times to create a set $V = \{v_{1}, v_{2}, \dots, v_{s + 1}\}$ and note that
$$
|H_{s + 1}| \geq |G| - (s + 1)(b - 1) \geq ((z - 1)\ell + 1)(b - 1)2k + z \geq z.
$$
This means that every vertex $v_{i} \in V$ has at least
$$
\frac{|H_{s + 1}|}{k} - ((z - 1)\ell + 1)(b - 1) \geq \frac{|H_{s + 1}|}{2k}
$$
edges to $H_{s + 1}$ in its own corresponding color.  Since $s = (z - 1)\ell$, by the pigeonhole principle, there exists a set of $z$ vertices $Z \subseteq V$, all with the same color, say red.  This means each vertex in $Z$ has at least $\frac{|H_{s + 1}|}{2k}$ red edges to $H_{s + 1}$.  By the definition of $z$, this subgraph contains the desired monochromatic (red) copy of $H$.
\end{proof}

Now suppose $H$ is not bipartite.  Let $t = R_{k}(H) - 1$ and let $m = m_{k}(H) = \min \{p :$ any $k$-coloring of $K_{\chi \cdot p}$ contains a monochromatic copy of $H \}$ where $K_{\chi \cdot p}$ is the complete multipartite graph with $\chi$ parts each of order $p$.

\begin{theorem}\label{Thm:MainNonBip}
Given positive integers $k \leq \ell$ and a graph $H$ with $\chi(H) = \chi \geq 3$ and $|H| = n$,
$$
(n - 1)(\chi - 1)^{\ell - 1} < ggr_{k, \ell}(H) \leq m\left( t2^{(\chi - 2) {\ell \choose k} + 1} + {\ell \choose k}(\chi - 1)m\right).
$$
\end{theorem}

\begin{proof}
For the lower bound, we iteratively construct a $k$-Gallai complete graph with no monochromatic copy of $H$ as follows.  Let $G_{1}$ be a graph on $n - 1$ vertices colored entirely with color $1$.  Given the graph $G_{i}$ construct the graph $G_{i + 1}$ by taking $\chi - 1$ disjoint copies of $G_{i}$ and inserting all edges between the copies in color $i + 1$.  Then $G_{\ell}$ certainly contains no monochromatic copy of $H$ and is $k$-Gallai for every $k \geq 1$.

For the upper bound, let $G$ be an $\ell$-coloring of a complete graph on 
$$
m\left( t2^{(\chi - 2) {\ell \choose k} + 1} + {\ell \choose k}(\chi - 1)m\right)
$$
vertices.  By assumption, there is a partition in which only $k$ colors appear on edges between the parts.  With $t = R_{k}(H) - 1$, there are at most $t$ parts in this partition.  Also, by the definition of $m$, there are at most $\chi - 1$ parts of the partition with order at least $m$.

Let $K$ be a $k$-subset of $[\ell]$.  Our goal is to construct, for each such subset $K$, $\chi - 1$ sets with $m$ vertices each, such that each $m$-set has only edges of colors in $K$ between each other and to a specified set of vertices.  Then if this specified set has at least $m$ vertices, the proof is complete by the definition of $m$.  Let $\{ K_{1}, K_{2}, \dots, K_{{\ell \choose k}} \}$ be the set of all these $k$-subsets of colors.  Define sets $T_{i, j}$ for $1 \leq i \leq {\ell \choose k}$ and $1 \leq j \leq \chi - 1$ to be empty sets.  We will fill these sets with vertices in the remainder of the proof.

Let $K_{i}$ be the set of colors used on edges between parts in our partition of $G$.  Let $A_{1}$ be a largest part of the partition and note that $|A_{1}| > m$.  For each part $A_{1}'$ of the partition, other than $A_{1}$, if $|A_{1}'| \geq m$, then place $m$ vertices from $A_{1}'$ into a set $T_{i, j}$ that is not full, with preference given to a set that is still empty.  The set $T_{i, j}$ then becomes \emph{full} since it contains $m$ vertices.  The other unused vertices of $A_{1}'$ are called \emph{wasted} since they are not used in any set $T_{i, j}$.  For smaller parts of the partition, those with order less than $m$, if possible, we collect several together with total order at least $m$ and add these to another set $T_{i, j}$ again that is not already full, with preference given to a set that is still empty.  Excess vertices here are again wasted.  The set $T_{i, j}$ then becomes full.  Finally, if the sum of the orders of the smaller parts is less than $m$, add vertices to a set $T_{i, j}$ that is not already full with preference given to a set that is not still empty.  Any excess vertices again get wasted.

By this process, we may be left with some sets $T_{i, j}$ that are full and at most one such set that is no longer empty but not yet full.  Also, fewer than $\frac{(\chi - 2)|G|}{\chi}$ vertices get wasted.  If there are fewer than $m$ vertices total outside $A_{1}$, say in a single part $A_{2}$, then at least one gets used to help fill a set $T_{i, j}$ so at most $m - 2$ vertices get wasted.

We then repeat this argument by considering the part $A_{1}$ as a $k$-Gallai complete graph and using the assumed partition.  If we let $K_{i'}$ be the set of colors used in this partition, we then use the vertices outside the largest part to fill sets $T_{i', j}$.

This process terminates when all the $\chi - 1$ sets $T_{i, j}$ get filled for some index $i$.  The process resulting in the most wasted vertices is when each pair of the iterated partitions of $n$ vertices (at the time) alternates first having $m - 1$ vertices outside the largest set, followed by a balanced bipartition.  Even this worst case results in fewer than $\frac{n}{2}$ wasted vertices per pair of iterations.  Then ${\ell \choose k}$ more iterations are performed in which $m - 1$ vertices are present outside the largest part.  Finally one balanced partition into $t$ parts would fill one of the final sets.  With $(\chi-2){\ell \choose k} + 1$ iterations, each roughly cutting the graph in half, the largest part in the final partition would have order at least
$$
\frac{m(t2^{(\chi - 2) {\ell \choose k} + 1} + {\ell \choose k}(\chi - 1)m)}{t2^{(\chi - 2) {\ell \choose k} + 1} + {\ell \choose k}(\chi - 1)m} = m.
$$
as desired.
\end{proof}


\end{document}